\documentclass[
reprint,showpacs,preprintnumbers,amsmath,amssymb,aps,pra]{revtex4-1}

\usepackage{graphicx,color}
\usepackage{dcolumn}
\usepackage{bm}
\usepackage{psfrag}
\usepackage{footnpag}
\usepackage{amsmath,amssymb}

\begin{document}

\preprint{arXiv/}

\title{Cyclic Operator Decomposition for Solving the Differential Equations}

\author{Ivan~Gonoskov}
\email{ivan.gonoskov@gmail.com}
\affiliation{
Institute of Applied Physics of the Russian Academy of Sciences\\
46 Ulyanov St., Nizhny Novgorod 603950, Russia}

\date{\today}

\begin{abstract}
	We present an approach how to obtain solutions of arbitrary linear operator equation for unknown functions. The particular solution can be represented by the infinite operator series (Cyclic Operator Decomposition), which acts the generating function. The method allows us to choose the cyclic operators and corresponding generating function selectively, depending on initial problem for analytical or numerical study. Our approach includes, as a particular case, the perturbation theory, but generally does not require inside any small parameters and unperturbed solutions. We demonstrate the applicability of the method to the analysis of several differential equations in mathematical physics, namely, classical oscillator, Schr\"odinger equation, and wave equation in dispersive medium.
\end{abstract}
%\pacs{02.30.Tb, 02.30.Hq, 02.30.Jr, 02.30.Lt, 02.30.Mv, 03.65.Db, 02.60.-x}

\maketitle

\section{Introduction}

	Various classical and quantum-mechanical problems in theoretical physics lead to the necessity of solving the linear operator equations for unknown functions and, in particular, the differential equations.  Exact non-trivial analytical solutions of these equations, which include finite combinations of elementary operations and special functions, are known only for a number of specific cases. However, there are many actual and important cases for which such exact solutions were not still obtained even by using severe approximations for the corresponding interaction operators. For the cases when exact solutions are unknown, some approximate methods are usually used. They can be conventionally divided into two types: (i) varieties of perturbation theory and (ii) numerical calculations (which generally are also based on perturbation theory). In spite of significant usefulness and applicability, these methods are not free from various limitations and disadvantages. The perturbation theory approaches may lead to divergent series, they need sometimes suitable unperturbed solutions, and, finally, they do not provide even estimations for precision in most cases (see \cite{Dyson-d}, \cite{Turbiner1984} and references therein). On the other hand, numerical schemes, which are from the very beginning approximate, usually also do not give reliable estimations for the precision (some reasonings can be found in \cite{FFM}). Moreover, they can hardly give an asymptotic behavior of the solutions at infinity.	Thus, the development of the general method which allows to overcome some of the above-mentioned difficulties is the main object of our study.
	
	In this manuscript we develop an approach based on the theory of Cyclic Operator Decomposition (COD), which gives the opportunities to obtain solutions (exact or approximate) of the differential equations with arbitrary operators. The particular solution can be represented by the infinite cyclic operators series, which acts the previously determined generating function. The cyclic operators and the corresponding generating function (COD components) can be specified through the given operators in the differential equation. Under the convergence requirement, these COD components can be chosen in different ways depending on the certain problem statement. The procedure differs from the using of Born series (or corresponding Neumann series) in the perturbation theory \cite{Corbett,Wellner,Perez} and S-matrix theory of Heisenberg, Feynmann and Dyson \cite{Dyson-S}. It can be understood easily by studying, for example, the difference between the formal definition of the generating function and Green's function (the last one is derived in some cases by using the operator resolvent formalism) \cite{Corbett,Perez,E-Reiss,Reiss-E,Mockel}. Generally, the proposed series does not require any small parameters or unperturbed solutions for the convergence. But, as a matter of fact, the procedure can be transformed, under some certain choice of COD components, to the "standard" perturbation theory with small parameters. For the potentials without strong singularities, with reasonable choice of the cyclic operators and generating function, the corresponding series usually has uniform convergence. Some additional features and advantages of our approach for analytical and numerical solving the differential equations are demonstrated in sections below. 
	
%	Below we present the theory of Cyclic Operator Decomposition. Then, We analyze convergence, possibility of finding the asymptotic behavior and estimation of accuracy by analyzing the equation of classical oscillator with time-dependent frequency. Next we demonstrate some particular choices of corresponding cyclic operators and generating functions for the problems, which are described by Schr\"odinger equation. Finally, we apply COD to the wave equation in dispersive medium.

\section{Theory of Cyclic Operator Decomposition}
	
	Let us start from the general case of operator equation for unknown function:
\begin{equation}\label{i1}
\hat{D}\psi=0.
\end{equation}
	Here $\hat{D}$ is an arbitrary given linear operator and $\psi$ is an unknown function, which can be a vector or matrix of arbitrary dimensionality. This equation can lead in particular cases to arbitrary linear differential equations, which are considered in examples below.
\newpage
	Let us consider a pair of operators $\hat{G}$ and $\hat{V}$, which are determined by the following condition:
\begin{equation}\label{i2}
\hat{D}=\hat{G}-\hat{V}.
\end{equation}
	
	Since the choice of this pair is partly optional, we impose additional conditions on the operator $\hat{G}$:
\begin{subequations}\label{i3}
\begin{align}
&\exists\hat{G}^{-1},\;\;\;\; &\text{that:}&\;\hat{G}\hat{G}^{-1}=\hat{I};\label{i-G-df}\\	
&\exists\psi_{g}\not\equiv{0},\;\;\;\; &\text{that:}&\;\hat{G}\psi_{g}\equiv{0},\label{i-f0-df}
\end{align}
\end{subequations}
where $\hat{I}$ is the identity operator. Any function $\psi_{g}$, which satisfies Eq.(\ref{i-f0-df}), will be called generating function. Now we can write the following equation:
\begin{equation}\label{i4}
\left(\hat{I}-\hat{G}^{-1}\hat{V}\right)\psi=\psi_{g}.
\end{equation}
As we can check, under the above-mentioned conditions for $\hat{G}$ and $\psi_{g}$, any solution of Eq.(\ref{i4}) fulfills Eq.(\ref{i1}). Equation (\ref{i4}) can be solved in terms of the following Cyclic Operator Decomposition:
\begin{equation}\label{i5}
\begin{split}	
\psi{}&=\left[\hat{I}+\hat{G}^{-1}\hat{V}+\hat{G}^{-1}\hat{V}\hat{G}^{-1}\hat{V}+...\right]\psi_{g}\\
&=\left[\hat{I}+\sum\limits_{n=1}^{\infty}\left(\hat{G}^{-1}\hat{V}\right)^{n}\right]\psi_{g}.
\end{split}
\end{equation}
	This is the exact particular solution of Eq.(\ref{i1}) with corresponding particular COD components determined by Eqs.(\ref{i2}) and (\ref{i3}). The solution makes sense only if the obtained series is convergent. This can be achieved in different cases depending on $\hat{G}^{-1}\hat{V}$ and $\psi_{g}$, for example, if we work in Banach space and corresponding operator norm is $\left\|\hat{G}^{-1}\hat{V}\right\|<1$. The convergence of some similar operator series was considered also in \cite{Corbett}.

	The theory can be easily generalized also to the case of the equations with given sources:
\begin{equation}\label{i7s}
\hat{D}\psi=\varphi,
\end{equation}	
where an arbitrary given function $\varphi$ describes the arbitrary sources. The unknown function $\psi$ could be found naturally if the inverse operator is known: $\psi=\hat{D}^{-1}\varphi$. However, the inverse operator $\hat{D}^{-1}$ can not be easily found for a number of problems. Then, we can write a solution of this equation analogously by using COD:
\begin{equation}\label{i7ss}
\psi{}=\left[\hat{I}+\sum\limits_{n=1}^{\infty}\left(\hat{G}^{-1}\hat{V}\right)^{n}\right]\left(\psi_{g}+\hat{G}^{-1}\varphi\right).
\end{equation}
In contrast to the case of Eq.(\ref{i1}), now we can choose $\psi_{g}\equiv{0}$ for some non-trivial particular solutions. Then we can obtain the particular solution of Eq.(\ref{i7s}), which corresponds to the following particular determination of the inverse operator in terms of COD:
\begin{equation}\label{i7sd}
\hat{D}^{-1}=\left[\hat{I}+\sum\limits_{n=1}^{\infty}\left(\hat{G}^{-1}\hat{V}\right)^{n}\right]\hat{G}^{-1}.
\end{equation}

	Important feature of the proposed theory is that, while the conditions Eqs.(\ref{i-G-df}) and (\ref{i-f0-df}) should be fulfilled and the convergence of the series is necessary, we still have a great freedom of choosing $\hat{G}$ and corresponding $\psi_{g}$. Generally, it gives us opportunities to obtain all the possible solutions of Eq.(\ref{i1}). Sometimes we can naturally choose $\hat{G}$ and $\psi_{g}$ in accordance, for example, with the corresponding initial conditions for Cauchy problem or boundary conditions for boundary-value problems. 
	
	In some cases, the exact solution Eq.(\ref{i5}) can be used naturally for obtaining the approximate solution with finite number of terms. It can be done, for example, when, starting from certain number $n$, the following conditions are satisfied: $\left\|\left(\hat{G}^{-1}\hat{V}\right)^{n}\psi_{g}\right\|\gg\left\|\left(\hat{G}^{-1}\hat{V}\right)^{n+1}\psi_{g}\right\|\gg...$. These are the sufficient conditions enabling one to derive the approximate solution with the prescribed accuracy. Further, the proposed method provides another advantage if one performs numerical calculations. According to the exact solution Eq.(\ref{i5}), we can use recurrent relations when calculating numerically the approximate solutions. In this case, the calculation of any next term in the corresponding series does not require more numerical resources than the calculation of the previous one.

\section{Examples}
	In this section we apply the proposed theory of Cyclic Operator Decomposition for the various cases of differential equations.
Let us first consider the Cauchy problem for the equation of classical oscillator. Note that this equation, if written in other variables, is the stationary one-dimensional Schr\"odinger equation with given energy, and it can be transformed also to the Riccati equation by using logarithmic substitution.
\begin{equation}\label{i8}
\begin{split}	
&\ddot{f}+\omega^{2}(t)f=0,\\
&{f}(t_{a})=a,\\
&\dot{f}(t_{b})=b,
\end{split}	
\end{equation}
where $f(t)$ is an unknown function, $\omega^{2}(t)$ is an arbitrary time-dependent frequency and $a,b$ are arbitrary constants. Here, it is natural to choose the components for COD as follows:

\begin{equation}\label{i9}
\begin{split}	
&\hat{G}=\frac{d^{2}}{dt^{2}},\;\;\;\;\;\;\;\psi_{g}=a+b(t-t_{a}),\\
&\hat{G}^{-1}=\int\limits_{t_{a}}^{t}d\tau_{1}\int\limits_{t_{b}}^{\tau_{1}}d\tau_{2},\;\;\;\; \hat{V}=-\omega^{2}(t).\\
\end{split}
\end{equation}
	
	If we fix (by our local convention, which we will use below) that we write for brevity the same variable upper limit of integration as the integration variable and determine the successive integration (step by step from right to left), we can write a simple expression for the solution: %(here we assume for simplicity that $a=1$, $b=0$):
\begin{widetext}
\begin{equation}\label{Schr-i13}
\begin{split}
f(t)=\left[a+b(t-t_{a})\right]-\int\limits_{t_{a}}^{t}dt\int\limits_{t_{b}}^{t}\omega^{2}(t)\left[a+b(t-t_{a})\right]dt+\int\limits_{t_{a}}^{t}dt\int\limits_{t_{b}}^{t}\omega^{2}(t)dt\int\limits_{t_{a}}^{t}dt\int\limits_{t_{b}}^{t}\omega^{2}(t)\left[a+b(t-t_{a})\right]dt-...
\end{split}
\end{equation}
\end{widetext} 
	It is important to note now, that the presented series (Eq.(\ref{Schr-i13})) has rapid uniform convergence at least in any interval, where $\omega^{2}(t)$ is bounded. For example, in the limited interval $[0,t]$ the rate of convergence for $f(t)$ can be estimated in the following way (we assume here for simplicity, that $t_{a}=t_{b}=0$, $b=0$, and the maximum of $\left|\omega^{2}(t)\right|$ in the corresponding interval is $C_{max}$):
\begin{equation}\label{Schr-i12}
(\text{n-th term of series})\leq\left|a\right|\frac{(C_{max})^{n}t^{2n}}{(2n)!}.
\end{equation}
In the same way, the convergence can be demonstrated for other different COD's, when the cyclic operators are bounded for the given generating functions in the given relevant interval. 
Moreover, if additionally $\omega(t)$ is a real function, $\omega^{2}(t)>\omega^{2}(0)\geq{0}$ and $t>0$, we have a decreasing alternating series for the above example, and we can estimate the precision of partial sum of the series by the value of the last term.

	Now we focus on some particular cases of $\omega^{2}(t)$. To demonstrate that it is possible to obtain a solution with any prescribed precision, we consider a case when \\$\omega^{2}(t)=[1-\frac{1}{2}\sin{t}]>0$ and $t_{a}=t_{b}=0$, $a=1$, $b=0$. Calculation of the first two terms in Eq.(\ref{Schr-i13}) gives
\begin{equation}\label{i14}
f(t)=\left[1-\frac{1}{2}(t^{2}-t+\sin{t})\right]+\delta.
\end{equation}
By calculating the third term in Eq.(\ref{Schr-i13}), we obtain $\delta<0.0273$ if we consider $t$ in the interval $[0,1]$.

	Sometimes we can find also the asymptotic behavior of the solution. As an example, we consider the case $\omega^{2}(t)=-t^{\alpha}$ and $t_{a}=t_{b}=0$, $a=1$, $b=0$,
where $\alpha$ is an arbitrary constant, $\alpha>-1$. Using again Eq.(\ref{Schr-i13}) we obtain exact solution in the following form:
\begin{equation}\label{i16}
\begin{split}
f(t)=1+&\frac{t^{\alpha+2}}{(\alpha+1)(\alpha+2)}\\
+&\frac{t^{2\alpha+4}}{(\alpha+1)(\alpha+2)(2\alpha+3)(2\alpha+4)}+...
\end{split}
\end{equation}
By analyzing the corresponding series we can obtain a simple upper estimate for the solution ($\forall{t}>0$):
\begin{equation}\label{i17}
f(t)<1+\frac{t^{\alpha+2}}{(\alpha+1)(\alpha+2)}\cdot\exp\left[\frac{2t^{\frac{1}{2}\alpha+1}}{\alpha+2}\right],
\end{equation}
which gives us the following asymptotic behavior at \\$t\rightarrow{\infty}$:
\begin{equation}\label{i18}
f(t)\propto\exp\left[\frac{2t^{\frac{1}{2}\alpha+1}}{\alpha+2}\right].
\end{equation} 
In this case, the same asymptotic can be found also from WKB theory (see, for example, quasiclassical approximation in \cite{LL-3,Schiff}).

	Let us now demonstrate the selective choice of cyclic operators. For that we consider stationary one-dimensional Schr\"odinger equation (we use below the units where $\hbar=1$, $m_{p}=1$):
\begin{equation}\label{i8sd}
\left[2E+\frac{d^{2}}{dx^{2}}-A\cdot{}e^{\beta{x}}\right]\psi(x)=0,
\end{equation}
where $\psi(x)$ is an unknown function, which describes quantum state with energy $E$ in the continuum; $A$, $\beta$ are arbitrary real constants. Without loss of generality (one can use scale transformations of Eq.(\ref{i8sd})) we can assume $2E=m^{2}, \beta=1$.
To find the solution of this equation, we can choose the components for COD in different ways. For example, if one interests in the behavior of $\psi(x)$ near $x=0$ and in small values of $E$, he can choose $\hat{G}=\frac{d^{2}}{dx^{2}}$ and use nearly the same technique as in Eqs.(\ref{i9}) and (\ref{Schr-i13}). However, this choice can be inconvenient for the analysis of the long-range behavior. 
Another variant of choosing the components for COD is the following (we use also Eq.(\ref{i7sd}) for the particular determination of $\hat{G}^{-1}$):
\begin{equation}\label{ii9}
\begin{split}	
&\hat{G}=m^{2}+\frac{d^{2}}{dx^{2}},\;\;\;\;\;\;\;\\
&\psi_{g}(x)=C_{1}e^{{i}mx}+C_{2}e^{-{i}mx},\\
&\hat{G}^{-1}=\left[\hat{1}+\sum\limits_{k=1}^{\infty}\left(\hat{G}_{0}^{-1}\hat{V}_{0}\right)^{k}\right]\hat{G}_{0}^{-1},\\
&\;\;\;\;\;\;\;\;\;\;\;\hat{G}_{0}^{-1}=\int\limits_{-\infty}^{x}dx\int\limits_{-\infty}^{x}dx,\\
&\;\;\;\;\;\;\;\;\;\;\;\hat{V}_{0}=-m^{2},\\
&\hat{V}=A\cdot{}e^{{x}},
\end{split}
\end{equation}
where $C_{1}$ and $C_{2}$ are arbitrary constants. To obtain the general solution, we consider the particular case of generating function $\tilde{\psi}_{g}(x)=e^{{i}mx}$ and corresponding particular solution $\psi_{p}(x)$. Then we derive by using a rule of infinite geometric series:
\begin{equation}\label{i9sd}
\begin{split}
\hat{G}^{-1}\hat{V}\tilde{\psi}_{g}(x)=\left[\hat{1}+\sum\limits_{k=1}^{\infty}\left(\hat{G}_{0}^{-1}\hat{V}_{0}\right)^{k}\right]\frac{A\cdot{}e^{(1+{i}m){x}}}{(1+{i}m)^{2}}\\
=\left[1+\sum\limits_{k=1}^{\infty}\left(\frac{m^{2}}{(1+{i}m)^{2}}\right)^{k}\right]\frac{A\cdot{}e^{(1+{i}m){x}}}{(1+{i}m)^{2}}=\frac{A\cdot{}e^{(1+{i}m){x}}}{(1+2{i}m)}.
\end{split}
\end{equation}
From here we obtain
\begin{equation}\label{i10sd}
\psi_{p}(x)=e^{{i}mx}\left[1+\sum\limits_{n=1}^{\infty}A^{n}e^{n{x}}\cdot{}\prod\limits_{k=1}^{n}\frac{1}{(k^{2}+2{i}m{k})}\right]
\end{equation}
and the general solution in the following form:
\begin{equation}\label{i10sd-g}
\psi(x)=C_{1}\cdot{}\psi_{p}(x)+C_{2}\cdot{}\psi^{*}_{p}(x)
\end{equation}
The corresponding series converges at any $A$ and real $m$. 

	In a similar way we can obtain solutions for multi-dimensional equations. Let us consider the stationary Schr\"odinger equation with potential surface $U(\vec{r})$, multi-dimensional Laplace operator $\Delta$, and energy $E$:
\begin{equation}\label{i26}
\left[\Delta+2(E-U(\vec{r}))\right]\psi(\vec{r})=0.
\end{equation} 	
Then we can choose $\hat{G}=\Delta$, $\hat{V}=-2(E-U(\vec{r}))$, and $\psi_{g}$ is any solution of $\Delta\psi_{g}=0$. From corresponding COD we can obtain a solution:
\begin{equation}\label{i27}
\psi(\vec{r})=\left[1+\Delta^{-1}V+\Delta^{-1}V\Delta^{-1}V+...\right]\psi_{g}=0.
\end{equation} 
%For the bounded potential $U(\vec{r})$
The inverse Laplace operator can be written, for example, as $\Delta^{-1}=-\hat{F}^{-1}k^{-2}\hat{F}$, where $\hat{F}$ is the Fourier transform operator and $k$ is an absolute value of the wave vector in this transform.\\ %This scheme can be realized recurrently for the calculation of wave function corrections.
Another choice of COD components can be better for the finding of the bound states with $E<0$. We can choose: $\hat{G}=2E+\Delta$, $\hat{V}=2U(\vec{r})$, $\psi_{g}$ is any solution of $[2E+\Delta]\psi_{g}=0$, and write the inverse operator $\hat{G}^{-1}$ for COD as follows:
\begin{equation}\label{i27-b}
\hat{G}^{-1}=\left(2E+\Delta\right)^{-1}=\hat{F}^{-1}\frac{1}{2E-k^{2}}\hat{F}.
\end{equation}
In this way, we can calculate in some cases the terms in the corresponding COD by evaluating the poles at imaginary values $k_{P}=\pm{i}\sqrt{-2E}$.

	Now we consider time-dependent three-dimensional Schr\"odinger equation to demonstrate other applications of the proposed method. Let us consider propagation of charged particle with arbitrary electromagnetic interactions (below $\vec{A}(\vec{r},t)$ is the vector potential, and $c=1$):
\begin{equation}\label{i19}
\left[i\frac{\partial}{\partial{t}}-\frac{1}{2}\left(i\nabla+\vec{A}(\vec{r},t)\right)^{2}-U(\vec{r},t)\right]\psi(\vec{r},t)=0.
\end{equation} 
Here, we can choose the components for COD in a variety of ways depending on peculiar properties of the interactions. One special choice is the following: %($t_{0}=0$):
\begin{equation}\label{i20}
\begin{split}	
&\hat{G}=i\frac{\partial}{\partial{t}},\;\;\;\;\;\;\;\psi_{g}=\psi(\vec{r},t_{0})=\psi_{0}(\vec{r}),\\
&\hat{G}^{-1}=-i\int\limits_{t_{0}}^{t}dt,\;\;\;\;\\
&\hat{V}=\frac{1}{2}\left(i\nabla+\vec{A}(\vec{r},t)\right)^{2}+U(\vec{r},t),\\
\end{split}
\end{equation}
where $\hat{V}$ corresponds to the time-dependent Hamiltonian. It gives the following solution:
\begin{equation}\label{i21}
\begin{split}	
\psi(\vec{r},t)=&\Bigl[1+(-i)\int\limits_{t_{0}}^{t}dt\hat{V}\\&+(-i)^{2}\int\limits_{t_{0}}^{t}dt\hat{V}\int\limits_{t_{0}}^{t}dt\hat{V}+...\Bigr]\psi_{0}(\vec{r}).
\end{split}
\end{equation} 

	This solution can be useful for the numerical calculations, namely, for the finding the propagator $\hat{P}$, which gives: $\psi(\vec{r},t+\Delta{t})=\hat{P}\psi(\vec{r},t)+O(\Delta{t}^{n})$. If we use additionally internal time-ordering, we can transform this expression to Dyson series (see \cite{Dyson-S,Dyson-S1}). 

	Finally, we consider the wave equation for electromagnetic waves in dispersive medium. We assume that an arbitrarily given operator $\hat{\epsilon}(\vec{r},t)$ which describes electric dispersion does not depend nonlinearly on the field, i.e. we still have linear problem. In this case the equation for unknown vector potential $\vec{A}(\vec{r},t)$ is the following (see for example \cite{LL-8}):  
\begin{equation}\label{i29}
\left[\frac{\partial}{\partial{t}}\left(\hat{\epsilon}(\vec{r},t)\frac{\partial}{\partial{t}}\right)-\nabla^{2}\right]\vec{A}(\vec{r},t)=0,
\end{equation} 	
with the initial conditions
\begin{equation}\label{i30}
\begin{split}	
&\vec{A}_{(t=0)}=\vec{S}(\vec{r}),\;\;\;\left(\frac{\partial\vec{A}}{\partial{t}}\right)_{(t=0)}=\vec{R}(\vec{r}).
\end{split}
\end{equation}	
We can choose the following components for COD:
\begin{equation}\label{i24t}
\begin{split}	
&\hat{G}=\frac{\partial}{\partial{t}}\left(\hat{\epsilon}(\vec{r},t)\frac{\partial}{\partial{t}}\right),\\
&\psi_{g}=\vec{S}(\vec{r})+\int\limits_{t_{0}}^{t}dt\;\hat{\epsilon}^{-1}(\vec{r},t)\vec{R}(\vec{r}),\\
&\hat{G}^{-1}=\int\limits_{t_{0}}^{t}dt\;\hat{\epsilon}^{-1}(\vec{r},t)\int\limits_{t_{0}}^{t}dt,\;\;\;\;\;\;\;\hat{V}=\nabla^{2}
\end{split}
\end{equation}	
Then, we can find the solution, which follows from the corresponding COD series:
\begin{equation}\label{i30}
\begin{split}	
\vec{A}(\vec{r},t)=\left[1+\sum\limits_{n=1}^{\infty}\left(\int\limits_{t_{0}}^{t}dt\;\hat{\epsilon}^{-1}(\vec{r},t)\int\limits_{t_{0}}^{t}dt\;\nabla^{2}\right)^{n}\right]\\
\times\left[\vec{S}(\vec{r})+\int\limits_{t_{0}}^{t}dt\;\hat{\epsilon}^{-1}(\vec{r},t)\vec{R}(\vec{r})\right].
\end{split}
\end{equation}	
For this solution the initial magnetic field is equal to $[\nabla\times\vec{S}(\vec{r})]$ and the initial electric field is equal to $\vec{R}(\vec{r})$.
\newpage

\section{Conclusions}

	In summary, we propose the theory of Cyclic Operator Decomposition, which allows one to obtain particular solutions of linear operator equations for unknown functions. In most cases it is possible to obtain all the possible solutions, which satisfy the given conditions. We demonstrate by some reasonings and particular examples that our approach has the following remarkable properties: (1) there is a freedom in choosing the COD components depending on the certain problem; (2) there is a rapid uniform convergence for most of the considered cases; (3) it is possible to find the asymptotic behavior of the solutions; (4) in many cases when one is analyzing the approximate solution, it is possible to estimate the accuracy; (5) the proposed approach gives good opportunities for efficient implementation of numerical calculations due to the recurrent relations that can be used in COD.

\section{ACKNOWLEDGMENTS}

Author would like to thank academician L. D. Faddeev, M. Yu. Emelin, M. Yu. Ryabikin, and A. A. Gonoskov for the useful and stimulating discussions.	

\bibliography{i-cod}

\end{document}